\documentclass[oneside,reqno,a4paper,12pt]{amsart}
\usepackage{amsfonts}
\usepackage{amsmath,amsfonts,amssymb,}
\usepackage[dvips]{epsfig}
\usepackage{CJK,color}
\usepackage{graphicx}

\newtheorem {theorem}{Theorem}[section]
\newtheorem {lemma}[theorem]{{\bf Lemma}}
\newtheorem {corollary}[theorem]{{\bf Corollary}}
\newtheorem {proposition}[theorem]{{\bf Proposition}}
\theoremstyle{remark}
\newtheorem {remark}{{\bf Remark}}[section]

\theoremstyle{problem}

\theoremstyle{definition}
\newtheorem {definition}{{\bf Definition}}[section]
\theoremstyle{plain} \numberwithin {equation}{section}

\makeindex

\begin{document}
\vspace{1cm}

\title[Singularity formation for Hall-MHD equations]{Singularity formation for the incompressible Hall-MHD equations without resistivity}
\author[Dongho Chae, \ \ Shangkun Weng]{\small Dongho Chae$^{\lowercase{1}}$, \ \ Shangkun Weng$^{\lowercase{2}}$}
%\thanks{$^*$ This research is supported by}
\thanks{$^1$ E-mail: dchae@cau.ac.kr. $^2$ E-mail: skwengmath@gmail.com. }
\maketitle
\begin{center} Department of Mathematics, Chung Ang University\\

          Seoul 156-756, Republic of Korea

\end{center}

\begin{abstract}
  In this paper we show that the incompressible Hall-MHD system without resistivity is not globally in time well-posed in any Sobolev space $H^{m}(\mathbb{R}^3)$ for any $m>\frac{7}{2}$. Namely, either the system is locally ill-posed in $H^{m}(\mathbb{R}^3)$, or it is locally well-posed, but  there exists an initial data in $H^{m}(\mathbb{R}^3)$, for which the  $H^{m}(\mathbb{R}^3)$ norm of solution blows-up in finite time if $m>7/2$.
  In the latter case we choose an axisymmetric initial data $u_0(x)=u_{0r}(r,z)e_r+ b_{0z}(r,z)e_z$ and $B_0(x)=b_{0\theta}(r,z)e_{\theta}$, and  reduce the system to the axisymmetric setting. If the convection term survives sufficiently long time, then the Hall term generates the singularity on the axis of symmetry and we have
  $ \lim\sup_{t\to t_*} \sup_{z\in \Bbb R} |\partial_z\partial_r b_\theta(r=0,z)|=\infty$ for some $t_*>0$.
\end{abstract}

\begin{center}
\begin{minipage}{5.5in}

\

\end{minipage}
\end{center}

\

\everymath{\displaystyle}
\newcommand {\eqdef }{\ensuremath {\stackrel {\mathrm {\Delta}}{=}}}

% math symbols
\def\Xint #1{\mathchoice
{\XXint \displaystyle \textstyle {#1}} %
{\XXint \textstyle \scriptstyle {#1}} %
{\XXint \scriptstyle \scriptscriptstyle {#1}} %
{\XXint \scriptscriptstyle \scriptscriptstyle {#1}} %
\!\int}
\def\XXint #1#2#3{{\setbox 0=\hbox {$#1{#2#3}{\int }$}
\vcenter {\hbox {$#2#3$}}\kern -.5\wd 0}}
\def\ddashint {\Xint =}
\def\dashint {\Xint -}
\def\clockint {\Xint \circlearrowright } % GOOD !
\def\counterint {\Xint \rotcirclearrowleft } % Good for Computer Modern !
\def\rotcirclearrowleft {\mathpalette {\RotLSymbol { -30}}\circlearrowleft }
\def\RotLSymbol #1#2#3{\rotatebox [ origin =c ]{#1}{$#2#3$}}

\def\be{\begin{eqnarray}}
\def\ee{\end{eqnarray}}
\def\p{\partial}
\def\no{\nonumber}
\def\e{\epsilon}
\def\de{\delta}
\def\De{\Delta}
\def\om{\omega}
\def\Om{\Omega}
\def\f{\frac}
\def\th{\theta}
\def\la{\lambda}
\def\b{\bigg}
\def\al{\alpha}
\def\La{\Lambda}
\def\ga{\gamma}
\def\Ga{\Gamma}
\def\ti{\tilde}
\def\Th{\Theta}
\def\si{\sigma}
\def\Si{\Sigma}
\def\bt{\begin{theorem}}
\def\et{\end{theorem}}
\def\bc{\begin{corollary}}
\def\ec{\end{corollary}}
\def\bl{\begin{lemma}}
\def\el{\end{lemma}}
\def\bp{\begin{proposition}}
\def\ep{\end{proposition}}
\def\br{\begin{remark}}
\def\er{\end{remark}}
\def\bd{\begin{definition}}
\def\ed{\end{definition}}
\def\bpf{\begin{proof}}
\def\epf{\end{proof}}

%\chapter{Subsonic irrotational flows in a two dimensional finitely long nozzle}\label{rectangular}\hspace*{\parindent}

\section{Introduction and main results}\label{introduction}

In this paper, we are concentrated on the singularity formation for the incompressible viscous Hall-MHD equations without resistivity. The incompressible viscous Hall-MHD equations without resistivity take the following form:
\be\label{hmhd}\begin{cases}
\p_t u+u\cdot\nabla u+\nabla p=(\nabla\times B)\times B +\nu\Delta u,\\
\text{div} u=0,\\
\p_t B-\nabla\times (u\times B)+\nabla\times ((\nabla\times B)\times B)=0,
\end{cases}
\ee
where $u(x,t)=(u_1(x,t), u_2(x,t), u_3(x,t))$ and $B(x,t)=(b_1(x,t), b_2(x,t),b_3(x,t))$, $(x,t)\in \Bbb R^3\times[0, \infty)$, are the fluid velocity and magnetic field. $\nu\geq 0$ is the viscosity, $\nu=0$ and $\nu>0$ correspond to the inviscid and viscous flow respectively. We will consider the Cauchy problem for (\ref{hmhd}), so we prescribe the initial data
\be\no
u(t=0,x)=u_0(x),\quad B(t=0,x)=B_0(x).
\ee
The initial data $u_0$ and $B_0$  satisfy the divergence free condition,
\be\no
\text{div}\, u_0(x)=\text{div} \,B_0(x)=0.
\ee
From the equations for the magnetic field $B$, it is easy to see that if one prescribes the divergence condition $\text{div}\, B_0=0$ on the initial data $B_0$, then $\text{div}\, B=0$ for later time.

Comparing with the well-known MHD system, the Hall term $\nabla\times ((\nabla\times B)\times B)$ is included due to the Ohm's law, which is believed to be a key issue for understanding magnetic reconnection. Note that the Hall term is quadratic in the magnetic field and involves the second order derivatives. Magnetic reconnection corresponds to a physical process in highly conducting plasmas in which the magnetic topology is rearranged and magnetic energy is converted to kinetic energy, thermal energy, and particle acceleration. During this process, the magnetic shear is large, the Hall term becomes dominant.  Lighthill \cite{lighthill} started the systematic study of the application of Hall-MHD on plasma, which is followed by \cite{campos}. One may refer to \cite{pm} for a physical review of the background for Hall-MHD.

There are many mathematical results on MHD system, for the existence of global weak solutions \cite{dl,st}, regularity criterion \cite{hx1,hx2} and global smooth small solutions \cite{lxz,xz}. The Hall-MHD has received little attention from mathematicians. The paper \cite{adfl} provided a derivation of Hall-MHD system from a two-fluids Euler-Maxwell system for electrons and ions, through a set of scaling limits. They also provided a kinetic formulation for the Hall-MHD,  and proved the existence of global weak solutions for the incompressible viscous resistive Hall-MHD system. The authors in \cite{cdl} obtained the local existence of smooth solutions for large data and global smooth solutions for small data to incompressible resistive, viscous or inviscid Hall-MHD model. Chae and Lee \cite{cl13} also established the blow-up criterion for classical solutions to the incompressible resistive Hall-MHD system. Contrary to the usual MHD, the global well-posedness question in the $2\f{1}{2}$ dimensional Hall-MHD is  still open. Note that $2\f{1}{2}$ dimensional Hall-MHD solution has been used in \cite{hg} to investigate the influence of the Hall term on the width of the magnetic islands of the tearing-mode. The temporal decay estimates for weak solutions to Hall-MHD system was established  by Chae and Schonbek \cite{cs}. They also obtained an algebraic decay rate for higher order Sobolev norms of solutions for small initial data.

In this paper we investigate the singularity formation for (\ref{hmhd}). Dreher, Ruban and Grauer \cite{drg} have discussed the possible spontaneous development of shock-type singularities in axisymmetric solutions of the ideal Hall-MHD system and performed numerical simulation to support their claim. In the following we rigorously prove that for the incompressible Hall-MHD system (\ref{hmhd}) without resistivity the solution cannot preserve initial data regularity in $H^m(\Bbb R^3)$, $m>7/2$. Either the solution breakdown the initial data regularity or uniqueness at the initial instant of moment, or if the solution survives  uniquely for a positive time, and if the convection term survives sufficiently long time, then a shock-type singularity in the magnetic field will develop in finite time, and this will also induce a singularity formation in the velocity field.

Now we start the mathematical setup of our problem and introduce the cylindrical coordinate
\be\no
r=\sqrt{x_1^2+x_2^2},\quad\theta=\arctan\f{x_2}{x_1},\quad z=x_3
\ee
and then investigate the axisymmetric solution to (\ref{hmhd}). In this case the velocity and magnetic field can be described as follows
\be\no
u(t,x)&=&u_r(r,z) e_r+ u_{\theta}(r,z)e_{\theta}+ u_z(r,z) e_z,\\\no
B(t,x)&=&b_r(r,z) e_r+ b_{\theta}(r,z)e_{\theta}+ b_z(r,z) e_z,\\\no
p(t,x)&=& p(t,r,z),
\ee
where
\be\no
e_r=(\cos\theta,\sin\theta,0),\quad e_{\theta}=(-\sin\theta,\cos\theta,0),\quad e_z=(0,0,1).
\ee

The Hall-MHD equation (\ref{hmhd}) can be written as the following equations in cylindrical coordinate
\be\no
&&\p_t u_r+\b((u_r\p_r+u_z\p_z) u_r-\f{u_{\theta}^2}{r}\b)+\p_r\b(p+\f{1}{2}(b_r^2+b_{\theta}^2+b_z^2)\b)\\\no
&&\quad\quad=\b((b_r\p_r+b_z\p_z) b_r-\f{b_{\theta}^2}{r}\b)+\nu(\p_r^2+\f{1}{r}\p_r+\p_z^2-\f{1}{r^2}) u_r,\\\no
&&\p_t u_{\theta}+\b((u_r\p_r+u_z\p_z) u_{\theta}+\f{u_r u_{\theta}}{r}\b)\\\no
&&\quad\quad=\b((b_r\p_r+b_z\p_z) b_{\theta}+\f{b_r b_{\theta}}{r}\b)+\nu(\p_r^2+\f{1}{r}\p_r+\p_z^2-\f{1}{r^2}) u_{\theta},\\\no
&&\p_t u_z+(u_r\p_r+u_z\p_z) u_z+\p_z\b(p+\f{1}{2}(b_r^2+b_{\theta}^2+b_z^2)\b)\\\no
&&\quad\quad =(b_r\p_r+b_z\p_z) b_z+\nu(\p_r^2+\f{1}{r}\p_r+\p_z^2) u_z,\\\no
&&\p_r u_r+\f{1}{r}u_r +\p_z u_z=0,\\\no
&&\p_t b_{r}+(u_r\p_r+u_z\p_z) b_r-(b_r\p_r+b_z\p_z)u_r-\f{\p}{\p z}(j_z b_r-j_r b_z)= 0,\\\no
&&\p_t b_{\theta}+\b((u_r\p_r+u_z\p_z) b_{\theta}+\f{b_r u_{\theta}}{r}\b)-\b((b_r\p_r+b_z\p_z)u_{\theta}+\f{u_r b_{\theta}}{r}\b)\\\no
&&\quad\quad+\b(\f{\p}{\p z}(j_{\theta}b_z-j_z b_{\theta})-\f{\p}{\p r}(j_r b_{\theta}-j_{\theta}b_r)\b)=0,\\\no
&&\p_t b_z+(u_r\p_r+u_z\p_z)b_z-(b_r\p_r+b_z\p_z)u_z+\f{1}{r}\f{\p}{\p r}\b(r(j_z b_r- j_r b_z)\b)=0,\\\no
&&\p_r b_r+\f{1}{r}b_r+\p_z b_z=0.
\ee
Here $j(t,x)=\nabla \times B=j_r(t,r,z) e_r+j_{\theta}(t,r,z) e_{\theta}+j_z(t,r,z)e_z$ and
\be\no
j_r=-\p_z b_{\theta},\quad j_{\theta}=\p_z b_r-\p_r b_z,\quad j_z=\f{1}{r}\p_r(r b_{\theta}).
\ee

From these equations, one can easily find that for any smooth solution $(u_r,u_{\theta}, u_z)$ and $(b_r,b_{\theta}, b_z)$, if initially one has
\be\label{initial}
u_{\theta}(0,r,z)=b_r(0,r,z)=b_z(0,r,z)=0,
\ee
then $u_{\theta}(t,r,z)=b_r(t,r,z)=b_{\theta}(t,r,z)\equiv 0$ for $t>0$. Hence $j_{\theta}\equiv 0$ and
\be\no
&\quad&\b(\f{\p}{\p z}(j_{\theta}b_z-j_z b_{\theta})-\f{\p}{\p r}(j_r b_{\theta}-j_{\theta}b_r)\b)=-\p_z (j_z b_{\theta})-\p_r(j_r b_{\theta})\\\no
&=&-(\p_r j_r+\p_z j_z) b_{\theta}-(j_r\p_r +j_z\p_z)b_{\theta}=-\f{2 b_{\theta}}{r} \p_z b_{\theta},
\ee
where we have used the fact that $\text{div} (\nabla\times B)=0$, so $\p_r j_r+\f{1}{r} j_r+\p_z j_z=0$. Finally under the initial condition (\ref{initial}) the above equations  reduce to
\be\label{a1}
\begin{cases}
\p_t u_r+(u_r\p_r+u_z\p_z) u_r+\p_r(p+\f{1}{2}b_{\theta}^2)=-\f{b_{\theta}^2}{r}+\nu(\p_r^2+\f{1}{r}\p_r+\p_z^2-\f{1}{r^2}) u_r,\\
\p_t u_z+(u_r\p_r+u_z\p_z) u_z+\p_z(p+\f{1}{2}b_{\theta}^2)=\nu(\p_r^2+\f{1}{r}\p_r+\p_z^2) u_z,\\
\p_r u_r+\f{1}{r}u_r +\p_z u_z=0,\\
\p_t b_{\theta}+(u_r\p_r+u_z\p_z) b_{\theta}-\f{u_r b_{\theta}}{r}-\f{2b_{\theta}}{r}\p_z b_{\theta} =0,\\
(u_r,u_z)(t=0,r,z)=(u_{0r},u_{0z})(r,z),\quad b_{\theta}(t=0,r,z)=b_{0\theta}(r,z).
\end{cases}
\ee

In this case, the vorticity $\omega(t,x)=\text{curl} u(x)=\omega_{\theta}(t,r,z) e_{\theta}=(\p_z u_r-\p_r u_z)(t,r,z) e_r$ satisfies the following equation
\be\no
&&\f{\p\omega_{\theta}}{\p t}+(u_r\p_r +u_z\p_z) \omega_{\theta}+2\f{b_{\theta}}{r}\p_z b_{\theta}-\f{u_r}{r} \omega_{\theta}=\nu(\p_r^2+\f{1}{r}\p_r+\p_z^2-\f{1}{r^2})\omega_{\theta},\\\no
&&\omega_{\theta}(t=0,r,z)=\omega_{0\theta}(r,z)=(\p_z u_{0r}-\p_r u_{0z})(r,z).
\ee
Define the new unknowns $\Omega=\f{\omega_{\theta}}{r}$ and $\Pi=\f{b_{\theta}}{r}$, then  one can check easily that $\Omega$ and $\Pi$ satisfy the following equations
\be\label{a101}
&&\f{\p \Omega}{\p t}+ (u_r \p_r+u_z\p_z) \Omega+2\Pi \p_z \Pi=\nu(\p_r^2+\f{3}{r}\p_r +\p_z^2) \Omega,\\\label{a102}
&&\f{\p \Pi}{\p t}+ (u_r\p_r+ u_z\p_z) \Pi-2\Pi\p_z \Pi=0,\\\label{a103}
&&\Omega(t=0,r,z)=\Omega_0(r,z):=\f{\omega_{0\theta}(r,z)}{r},\\\label{a104}
&&\Pi(t=0,r,z)=\Pi_0(r,z):=\f{b_{0\theta}(r,z)}{r}.
\ee

%First we establish the local in time existence of smooth solution to (\ref{a1}) with smooth initial data.
%\bt\label{existence}
%{\it
%For any smooth initial data $u_0(x)=u_{0r}(r,z)e_r+ u_{0z}(r,z) e_z\in H^4(\mathbb{R}^3)$ and $B_0(x)=b_{0\theta}(r,z) e_{\theta}\in H^4(\mathbb{R}^3)$, which also satisfies $\f{b_{0\theta}}{r}\in L^{\infty}([0,\infty)\times \mathbb{R})$ and some compatibility conditions
%\be\no
%u_{0r}(0,z)=b_{0\theta}(0,z)=0,\\\no
%\p_r^2 u_{0r}(0,z)=\p_r^2 b_{0\theta}(0,z)=0,
%\ee
%then there exists a smooth solution $(u(t,x), B(t,x))\in L^{\infty}([0,T_0]; H^4(\mathbb{R}^3))\cap Lip([0,T_0]; H^2(\mathbb{R}^3))$, which has the following form
%\be\no
%u(t,x)= u_r(t,r,z)e_r+ u_z(t,r,z) e_z,\quad B(t,x)=b_{\theta}(t,r,z)e_{\theta},
%\ee
%where $T_0=T_0(\|u_0\|_{H^4(\mathbb{R}^3)}, \|B_0\|_{H^4(\mathbb{R}^3)})$. Indeed, $(u_r,u_z, b_{0\theta})(t,r,z)$ is a smooth solution to (\ref{a1}).
%}
%\et

%\br\label{existenceforgeneraldata}
%{\it
%The existence of a local in time smooth solution for (\ref{hmhd}) for general initial data is not clear so far. The Hall term makes it difficult to get a close energy estimate.
%}
%\er

We refer two closely related results on the axisymmetric solution to the usual MHD or Hall-MHD system. Lei \cite{l} has showed that the existence of global in time smooth solution to the incompressible viscous MHD without resistivity for some special axisymmetric data $u_0=u_{0r} e_r+ u_{0z} e_z$ and $B_0= b_{0\theta} e_{\theta}$. The result in \cite{fhn} established the existence of global smooth solution to the incompressible viscous, resistive Hall-MHD system with same initial data. Our first main result is the formation singularity for the incompressible viscous Hall-MHD without resistivity.

\bt\label{a}({\it Viscous case $\nu=1$}).
{\it
The incompressible viscous Hall-MHD system without resistivity (\ref{hmhd}) is not globally well-posedness in any Sobolev space $H^{m}(\mathbb{R}^3)$ for $m>\frac{7}{2}$. There exists smooth initial data $u_0(x)=u_{0r}(r,z) e_{r}+ u_{0z}(r,z) e_z\in C_c^{\infty}(\mathbb{R}^3)$, $B_0(x)=b_{0\theta}(r,z) e_{\theta}\in C_c^{\infty}(\mathbb{R}^3)$ with $\Pi_0(r,z)\in L^{\infty}(\mathbb{R}^3)$ such that if there is a local in time smooth solution $(u,B)(t,x)$ to (\ref{hmhd}) with initial data $(u_0,B_0)$, then $(u,B)$ must blow up in finite time. Indeed, one can choose $(u_0,B_0)$ such that $y_0:=\p_z\Pi_0(0,0)=\p_{rz}^2 b_{0\theta}(0,0)\geq 10^4 C_*^2, t_0=\f{4}{y_0}$ and $J_0:=\Pi_0(0,0)>0$, where $C_*$ depends only on $\|u_0\|_{H^2(\mathbb{R}^3)}, \|B_0\|_{H^1(\mathbb{R}^3)}$ and $\|\Pi_0\|_{L^{\infty}(\mathbb{R}^3)}$, then
\be\no
\lim\sup_{t\to t_0}\sup_{z\in \Bbb R} |\p_z\Pi(t,0,z)| =\infty.
\ee
Moreover, the velocity field also blows up
\be\no
\lim\sup_{t\rightarrow t_0} \sup_{z\in \Bbb R}\left|\left(\p_t \Omega+(u_r\p_r+u_z\p_z)\Omega-\nu(\p_r^2+\f{3}{r}\p_r+\p_z^2)\Omega\right)(t,0,z)\right|=\infty.
\ee
}
\et

\br\label{remark1}
{\it As will be shown in the proof below, the singularity occurs on the axis if the local well-posednes is done. The blow-up happens on the second order derivative of $b_{\theta}$ and the third derivative of the velocity field. Whether the solution can blow-up off the axis is not clear yet.
}
\er
\br\label{remark11}
{\it Due to the Hall term it seems difficult to show that the local in time existence of smooth solution to (\ref{hmhd}).  We could not rule out the possibility at this moment that (\ref{hmhd}) is locally ill-posed(see Remark 3.1 of \cite{cdl}).}
\er

\br\label{remark2}
{\it
If one consider the equations (\ref{hmhd}) with only partial viscosity in the $z$-direction, i.e. replace $\Delta u$ by $\p_z^2 u$, then Theorem \ref{a} is still true. We will indicate the corresponding modification in the following section.
}
\er

The second result concentrates on the singularity formation for the inviscid Hall-MHD system without resistivity.
\bt\label{b} ({\it Invisicd case $\nu=0$}).
{\it
The incompressible inviscid Hall-MHD system without resistivity (\ref{hmhd}) is not globally well-posedness in any Sobolev space $H^{m}(\mathbb{R}^3)$ for $m>\frac{7}{2}$. There exists smooth initial data $u_0(x)=u_{0r}(r,z) e_{r}+ u_{0z}(r,z) e_z\in C_c^{\infty}(\mathbb{R}^3)$, $B_0(x)=b_{0\theta}(r,z) e_{\theta}\in C_c^{\infty}(\mathbb{R}^3)$ with $(\Omega_0,\Pi_0)(r,z)\in L^{1}(\mathbb{R}^3)\cap L^{\infty}(\mathbb{R}^3)$ such that if there is a local in time smooth solution $(u,B)(t,x)$ to (\ref{hmhd}) with initial data $(u_0,B_0)$, then $(u,B)$ must blow up in finite time. Indeed, one can choose $(u_0,B_0)$ such that $y_0:=\p_z\Pi_0(0,0)=\p_{rz}^2 b_{0\theta}(0,0)\geq 4 \tilde{C}_{*}^{1/2}, t_{\sharp}=\f{4}{y_0}$ and $J_0:=\Pi_0(0,0)>0$, where $\tilde{C}_{*}$ depends only on $\|\Omega_0\|_{L^1\cap L^{\infty}}+\|\Pi_0\|_{L^1\cap L^{\infty}}$, where $\|f\|_{L^1\cap L^{\infty}}:=\|f\|_{L^1(\mathbb{R}^3)}+\|f\|_{L^{\infty}(\mathbb{R}^3)}$, then
\be\no
\lim\sup_{t\to t_{\sharp}}\sup_{z\in \Bbb R} |\p_z\Pi(t,0,z)| =\infty.
\ee
Moreover, the velocity field also blows up
\be\no
\lim\sup_{t\rightarrow t_{\sharp}} \sup_{z\in \Bbb R}\left|\left(\p_t \Omega+(u_r\p_r+u_z\p_z)\Omega\right)(t,0,z)\right|=\infty.
\ee
}
\et

The paper will proceed as follows. In section \ref{apriori}, we will give some a priori estimates on the smooth solutions to (\ref{a1}). Then we prove Theorem \ref{a} and \ref{b} in the last section.

\section{Some a priori estimates for solutions to (\ref{a1})}\label{apriori}

\subsection{A priori estimates: viscous case $\nu=1$}
First we give some a priori estimates for solutions to (\ref{a1}). The following lemma shows that the maximum principle for $\Pi$. The proof is easy, we omit the details.
\bl\label{mp}
{\it
For any smooth solution $(u_r,u_z,b_{\theta},p)$ to (\ref{a1}) with initial data $u_0(x)=u_{0r}(r,z) e_r+ u_{0z}(r,z) e_z\in C_c^{\infty}(\mathbb{R}^3)$, $B_0(x)=b_{0\theta}(r,z) e_{\theta}\in C_c^{\infty}(\mathbb{R}^3)$ satisfying $\Pi_0(r,z)\in L^{\infty}(\mathbb{R}^3)$, then we have
\be\no
\|\Pi(t,r,z)\|_{L^{\infty}}\leq \|\Pi_0(r,z)\|_{L^{\infty}}.
\ee
If $\Pi_0\in L^2(\mathbb{R}^3)$, then
\be\no
\|\Pi(t,\cdot)\|_{L^2(\mathbb{R}^3)}=\|\Pi_0\|_{L^2(\mathbb{R}^3)}.
\ee
}
\el

\bl\label{enstrophy}{\it ($L^2$ estimate of $\Omega$.) Assume that the initial data $(u_0,B_0)$ satisfy $u_0\in H^2(\mathbb{R}^3),B_0\in H^1(\mathbb{R}^3)$ and $\Pi_0\in L^{\infty}$. Then we have the following estimate for $\Omega$
\be\no
&\quad&\|\Omega(t,\cdot)\|_{L^2}^2+\int_0^t \|\nabla \Omega(s,\cdot)\|_{L^2}^2 ds+2\pi\int_0^t\int_{\Bbb R} |\Omega(s,0,z)|^2 dzds\\\no
&\leq& C_1(\|u_0\|_{H^2(\mathbb{R}^3)}, \|B_0\|_{H^1(\mathbb{R}^3)}, \|\Pi_0\|_{L^{\infty}})(1+t).
\ee
}
\el

\bpf
By (\ref{a101}), one can easily obtain the $L^2$ estimate for $\Omega$
\be\no
&\quad&\f{d}{dt} \|\Omega\|_{L^2}^2 + \|\nabla \Omega\|_{L^2}^2+2\pi\int_{\mathbb{R}}|\Omega(t,0,z)|^2 dz\\\no
&=& -\int_{\Bbb R^3} \Omega \p_z \Pi^2 dx =\int_{\Bbb R^3} \Pi^2 \p_z\Omega dx\\\no
&\leq& \|\Pi\|_{L^{\infty}} \|\Pi\|_{L^2} \|\p_z \Omega\|_{L^2}\leq 4\|\Pi\|_{L^{\infty}}^2 \|\Pi\|_{L^2}^2+\f{1}{2}\|\p_z\Omega\|_{L^2}^2.
\ee
Hence we obtain
\be\no
&\quad&\f{d}{dt} \|\Omega\|_{L^2}^2+\|\nabla \Omega\|_{L^2}^2+2\pi \int_{\mathbb{R}}|\Omega(t,0,z)|^2 dz\\\no
&\leq& 4\|\Pi\|_{L^{\infty}}^2 \|\Pi\|_{L^2}^2\leq 4 \|\Pi_0\|_{L^{\infty}}^2 \|\Pi_0\|_{L^2}^2.
\ee
This will imply the following estimate for $\Omega$
\be\no
&\quad&\|\Omega(t,\cdot)\|_{L^2}^2+\int_0^t \|\nabla\Omega(s,\cdot)\|_{L^2}^2 ds +2\pi\int_0^t \int_{\Bbb R}|\Omega(t,0,z)|^2 dz ds \\\no
&\leq& \|\Omega_0\|_{L^2}^2+ 4 \|\Pi_0\|_{L^{\infty}}^2 \|\Pi_0\|_{L^2}^2 t\\\no
&\leq& \|u_0\|_{H^2}^2+ 4 \|\Pi_0\|_{L^{\infty}}^2\|B_0\|_{H^1}^2 t\leq C_1(\|u_0\|_{H^2},\|\Pi\|_{L^{\infty}},\|B_0\|_{H^1})(1+t).
\ee

\epf
\br\label{remark3}
{\it If one consider the equations (\ref{hmhd}) with only partial viscosity in the $z$-direction, i.e. replace $\Delta u$ by $\p_z^2 u$, then we still have the following estimate
\be\label{remark31}
&\quad&\|\Omega(t,\cdot)\|_{L^2}^2+\int_0^t \|\p_z \Omega(s,\cdot)\|_{L^2}^2 ds\\\no
&\leq& C_1(\|u_0\|_{H^2(\mathbb{R}^3)}, \|B_0\|_{H^1(\mathbb{R}^3)}, \|\Pi_0\|_{L^{\infty}})(1+t).
\ee
}
\er

We also need the following estimate for $\f{u_r}{r}$. This estimate has been appeared in Lemma 3.1 in \cite{l}.
\bl\label{ur}
{\it
The following estimate holds for $\f{u_r}{r}$:
\be\label{ur1}
\int_0^t\left\|\f{u_r}{r}(s,\cdot)\right\|_{L^{\infty}}^2 ds &\leq& \sup_{0\leq s\leq t}\|\Omega(s,\cdot)\|_{L^2} \int_0^t \|\p_z \Omega(s,\cdot)\|_{L^2}^2 ds\leq  C_* (1+t)^{3/2} t^{1/2},
\ee
where $C_*$ depends only on $\|u_0\|_{H^2(\mathbb{R}^3)}, \|B_0\|_{H^1(\mathbb{R}^3)}, \|\Pi_0\|_{L^{\infty}}$.

}
\el

\bpf
For the convenience of the reader we give a sketch of proof. For more details of the proof, one may refer to \cite{l}. By the divergence free condition, $\p_r (r u_r)+\p_z (r u_z)=0$, one can introduce a stream function $\psi^{\theta}$ such that
\be\no
u_r =-\p_z \psi_{\theta},\quad u_z=\f{1}{r}\p_r(r\psi_{\theta}).
\ee
Since $\omega_{\theta}=\p_z u_r-\p_r u_z$, we have
\be\no
-(\p_r^2 +\f{1}{r}\p_r+\p_z^2-\f{1}{r^2})\psi_{\theta}= \omega_{\theta}.
\ee
Setting $\varphi=\f{\psi_{\theta}}{r}$, then it is easy to see that
\be\no
-(\p_r^2+\f{3}{r}\p_r+\p_z^2)\varphi=\Omega.
\ee
As in \cite{l}, the second order operator $(\p_r^2+\f{3}{r}+\p_z^2)$ can be interpreted as the Laplace operator in 5-dimensional space. We introduce
\be\no
y=(y_1,y_2,y_3,y_4, z),\quad r=\sqrt{y_1^2+y_2^2+y_3^2+y_4^2},\quad \Delta_y=(\p_r^2+\f{3}{r}\p_r+\p_z^2).
\ee
Hence we have $\varphi=(-\Delta_y)^{-1} \Omega$. To get an estimate of $\left\|\f{u_r}{r}\right\|_{L^{\infty}}$, by a simple interpolation inequality $\|f\|_{L^{\infty}}^2\leq \|\nabla f\|_{L^2}\|\nabla^2 f\|_{L^2}$, we have
\be\no
\int_0^t \left\|\f{u_r}{r}(s,\cdot)\right\|_{L^{\infty}}^2 ds&=& \int_0^t \|\p_z \varphi(s,\cdot)\|_{L^{\infty}}^2 ds\\\no
&\leq& \int_0^t \|\nabla \p_z \varphi(s,\cdot)\|_{L^2}\|\nabla^2 \p_z \varphi(s,\cdot)\|_{L^2} ds.
\ee
By simple calculations, one has
\be\no
|\nabla_y^2 \varphi|^2\simeq |\p_r^2 \varphi|^2+|\f{1}{r}\p_r \varphi|^2+|\p_z^2 \varphi|^2+|\p_{rz}^2 \varphi|^2
\ee
and
\be\no
\int |\nabla^2 \varphi|^2 dx&\leq& C\int_{-\infty}^{\infty} \int_0^{\infty} \b(|\p_r^2 \varphi|^2+|\f{1}{r}\p_r \varphi|^2+|\p_z^2 \varphi|^2+|\p_{rz}^2 \varphi|^2\b) r dr dz\\\no
&=&C\int_{-\infty}^{\infty} \int_0^{\infty} \b(|\p_r^2 \varphi|^2+|\f{1}{r}\p_r \varphi|^2+|\p_z^2 \varphi|^2+|\p_{rz}^2 \varphi|^2\b)w(r) r^3 dr dz\\\no
&\leq& C\int_{-\infty}^{\infty} \int_0^{\infty} |\nabla_y^2 \varphi|^2w(r) r^3 dr dz= \int_{-\infty}^{\infty} \int_0^{\infty} |\nabla_y^2 (-\Delta_y)^{-1}\Omega|^2w(r) r^3 dr dz\\\no
&=& C\int|\nabla_y^2 (-\Delta_y)^{-1}\Omega|^2 w(r) dy\\\no
&\leq& C\int |\Omega|^2 w(r) dy= \int |\Omega|^2 dx,
\ee
where $w(r)=r^{-2}$ and in the last step we have used Lemma 2 in \cite{hll}. See also Corollary 2 in \cite{cl02} for a similar weighted estimate for a singular integral operator.

Similarly, we also have
\be\no
\int |\nabla^2 \p_z\varphi|^2 dx\leq \int |\p_z \Omega|^2 dx.
\ee
Hence
\be\no
\int_0^t \left\|\f{u_r}{r}(s,\cdot)\right\|_{L^{\infty}}^2 ds &\leq & \int_0^t \|\nabla \p_z \varphi(s,\cdot)\|_{L^2}\|\nabla^2 \p_z \varphi(s,\cdot)\|_{L^2} ds\\\no
&\leq& \int_0^t \|\Omega(s,\cdot)\|_{L^2} \|\p_z \Omega(s,\cdot)\|_{L^2} ds\\\no
&\leq&C(1+t) \b(\int_0^t \|\p_z \Omega(s,\cdot)\|_{L^2}^2 ds\b)^{1/2} t^{1/2}\leq C_*(1+t)^{3/2} t^{1/2}.
\ee

\epf
\br\label{remark4}
{\it If one consider the equations (\ref{hmhd}) with only partial viscosity in the $z$-direction, i.e. replace $\Delta u$ by $\p_z^2 u$, then we still have the following estimate
\be\label{remark41}
\int_0^t\left\|\f{u_r}{r}(s,\cdot)\right\|_{L^{\infty}}^2 ds &\leq& \sup_{0\leq s\leq t}\|\Omega(s,\cdot)\|_{L^2} \int_0^t \|\p_z \Omega(s,\cdot)\|_{L^2}^2 ds\leq  C_* (1+t)^{3/2} t^{1/2}.
\ee
}
\er

\subsection{A priori estimates: Inviscid case $\nu=0$.}
In this case, then the equations satisfied by $\Omega$ and $\Pi$ will reduce to
\be\label{b1}
\begin{cases}
\p_t \Omega+ (u_r\p_r +u_z\p_z)\Omega+ 2\Pi\p_z \Pi=0,\\
\p_t \Pi +(u_r\p_r +u_z\p_z) \Pi-2\Pi\p_z\Pi=0,\\
(\Omega,\Pi)(t=0,r,z)=(\Omega_0,\Pi_0)(r,z).
\end{cases}
\ee
Putting $\Gamma=\Omega+\Pi$, it is easy to see that
\be\label{b2}
\begin{cases}
\p_t \Gamma+(u_r\p_r+u_z\p_z) \Gamma=0,\\
\Gamma(t=0,r,z)=\Omega_0(r,z)+\Pi_0(r,z):=\Gamma_0(r,z).
\end{cases}
\ee
This simple, but important observation plays a key role in our following argument. Note that (\ref{b2}) indeed comes from (\ref{hmhd}) with $\nu=0$ by observing that $R=\text{curl} \, u+ B$ satisfies the following equation
\be\label{b3}
\p_t R+ u\cdot \nabla R-R\cdot\nabla u=0.
\ee

\bl\label{b4}
{\it
For any smooth solution $(u_r,u_z,b_{\theta},p)$ to (\ref{a1}) with initial data $u_0(x)=u_{0r}(r,z) e_r+ u_{0z}(r,z) e_z\in C_c^{\infty}(\mathbb{R}^3)$, $B_0(x)=b_{0\theta}(r,z) e_{\theta}\in C_c^{\infty}(\mathbb{R}^3)$ satisfying $(\Omega_0,\Pi_0)(r,z)\in L^1\cap L^{\infty}$, then we have
\be\label{b5}
\|\Pi(t,r,z)\|_{L^1\cap L^{\infty}}&\leq& \|\Pi_0(r,z)\|_{L^1\cap L^{\infty}},\\\label{b6}
\|\Omega(t,r,z)\|_{L^1\cap L^{\infty}}&\leq& \|\Omega_0(r,z)\|_{L^1\cap L^{\infty}}+\|\Pi_0(r,z)\|_{L^1\cap L^{\infty}}.
\ee
}
\el

Next we need the following inequality, which comes from the Biot-Savart law and has been proved in \cite{sy} long time ago. One can refer to Lemma 2 in \cite{danchin} for more details.
\bl\label{b6}
{\it
There exists a universal constant $C_2$ such that
\be\label{b7}
|u_r(t,x)|\leq C_2 \int_{\mathbb{R}^3} \min\b(1,\f{r}{|x'-x|}\b)\f{|\omega_{\theta}(t,x')|}{|x-x'|^2} dx',
\ee
which yields
\be\label{b8}
\f{|u_r(t,x)|}{r} \leq 2C_2 \int_{\mathbb{R}^3} \f{1}{|x-x'|^2} \f{|\omega_{\theta}(x')|}{r'}d x'.
\ee
Note that here we use the notation $u_r(t,x):=u_r(t,\sqrt{x_1^2+x_2^2}, x_z)$.
}
\el

From (\ref{b8}), we have for any $t>0$
\be\no
\left\|\f{u_r}{r}(t,\cdot)\right\|_{L^{\infty}(\mathbb{R}^3)}&\leq& 2C_2 \int_{\mathbb{R}^3} \f{1}{|x-x'|^2} |\Omega(x')| dx'\\\no
&= & 2 C_2 \b(\int_{ |x-x'|\leq 1}+ \int_{|x-x'|> 1}\b)\f{1}{|x-x'|^2} |\Omega(x')| dx'\\\no
&\leq& 2 C_2 \b(\|\Omega\|_{L^{\infty}}\int_{|x-x'|\leq 1}\f{1}{|x-x'|^2} d x'+ \|\Omega\|_{L^1(\mathbb{R}^3)}\b)\\\no
&\leq& C_3 \|\Omega\|_{L^1\cap L^{\infty}}\\\label{b9}
&\leq& C_3\b(\|\Omega_0(r,z)\|_{L^1\cap L^{\infty}}+\|\Pi_0(r,z)\|_{L^1\cap L^{\infty}}\b):=\tilde{C}_{*},
\ee
where $C_3$ is also a universal constant.

\section{Singularity formation}

\subsection{Viscous case $\nu=1$.}

{\it Proof of Theorem \ref{a}.} Suppose the incompressible viscous Hall-MHD system without resistivity (\ref{hmhd}) is globally well-posedenss in any Sobolev space $H^{m}(\mathbb{R}^3)$ for $m> \f{7}{2}$.  We will derive contradiction to this. \\
If the system is locally ill-posed, then we are done, and nothing to prove.  Therefore, we assume that (\ref{hmhd})
is locally in time well-posed in $H^{m}(\mathbb{R}^3)$ for $m> \f{7}{2}$. Namely there exists $T>0$ such that a unique solution $(u, B)\in \{C([0, T);H^{m}(\mathbb{R}^3)\}^2$ exists.
In the following, we will choose a special class of smooth axisymmetric initial data with the form $u_0(x)= u_{0r}(r,z)e_{r}+u_{0z}(r,z) e_z\in ( C_c^{\infty}(\mathbb{R}^3))^3$ and $B_0(x)=b_{0\theta}(r,z) e_{\theta} \in ( C_c^{\infty}(\mathbb{R}^3))^3$, such that the corresponding solution $(u,B)(t,x)$ to (\ref{hmhd}) will develop in finite time  a singularity for the magnetic field, which will also induce a singularity in the velocity field. Hence we can conclude that the Hall-MHD system (\ref{hmhd}) is not global well-posedness in any Sobolev space $H^{m}(\mathbb{R}^3)$ for $m>\f{7}{2}$.

For Hall-MHD system with initial data $u_0(x)= u_{0r}(r,z)e_{r}+u_{0z}(r,z) e_z$ and $B_0(x)=b_{0\theta}(r,z) e_{\theta}$, by uniqueness, we can show that the corresponding solution $(u,B)(t,x)$ should be axisymmetric and has the form
\be\no
u(t,x)=u_{r}(t,r,z) e_r+ u_{z}(t,r,z) e_{z},\quad B(t,x)= b_{\theta}(t,r,z) e_{\theta},
\ee
where $(u_r, u_z, b_{\theta})(t,r,z)$ should solve the system (\ref{a1}) with initial data $(u_{0r}, u_{0z}, b_{0\theta})$. Indeed, for any $\alpha \in [0, 2\pi)$, we define the following change of coordinate
\be\no
\left(\begin{array}{l}
y_1\\
y_2\\
y_3
\end{array}\right):= A \left(\begin{array}{l}
x_1\\
x_2\\
x_3
\end{array}\right)=\left(\begin{array}{lll}
\cos \alpha & \sin\alpha & 0\\
-\sin\alpha & \cos\alpha & 0\\
0 & 0& 1
\end{array}\right) \left(\begin{array}{l}
x_1\\
x_2\\
x_3
\end{array}\right).
\ee
Setting
\be\no
\tilde{u}(t,y)= A u(t, A^{-1}y),\quad \tilde{B}(t,y)=A u(t, A^{-1}y),\quad \tilde{p}(y)= \tilde{p}(A^{-1}y),
\ee
then it is easy to verify that $(\tilde{u}(t,y),\tilde{B}(t,y), \tilde{p}(t,y))$ solves (\ref{hmhd}) with initial data
\be\no
\tilde{u}(t=0,y)= A u_0(A^{-1}y),\quad \tilde{B}(t=0,y)=A B_0(A^{-1}y).
\ee
By the axisymmetric property of $(u_0(x), B_0(x))$, we have $A u_0(A^{-1}y)= u_0(y), A B_0(A^{-1}y)=B_0(y)$. Hence by uniqueness of (\ref{hmhd}), we have
\be\no
\tilde{u}(t,y)\equiv u(t,y),\quad \tilde{B}(t,y)\equiv B(t,y),\quad \tilde{p}(t,y)\equiv p(t,y).
\ee
Since $\alpha\in [0,2\pi)$ is arbitrary, we find that $(u,B,p)(t,x)$ must be axisymmetric and  is of the form
\be\no
u(t,x)=u_r(t,r,z)e_r+u_z(t,r,z)e_z,\quad B(t,x)=b_{\theta}(t,r,z) e_{\theta},\quad p(t,x)=p(t,r,z)
\ee
where $(u_r,u_z, b_{\theta},p)(t,r,z)$ solves the problem (\ref{a1})(see lines below (\ref{initial})).

Hence the a priori estimates established in section \ref{apriori} hold for $(u_r,u_z,b_{\theta})(t,r,z)$. In particular, we have the following estimate
\be\no
\int_0^t\left\|\f{u_r}{r}(s,\cdot)\right\|_{L^{\infty}}^2 ds &\leq& \sup_{0\leq s\leq t}\|\Omega(s,\cdot)\|_{L^2} \int_0^t \|\p_z \Omega(s,\cdot)\|_{L^2}^2 ds\leq  C_* (1+t)^{3/2} t^{1/2},
\ee
where $C_*$ depends only on $\|u_0\|_{H^2(\mathbb{R}^3)}, \|B_0\|_{H^1(\mathbb{R}^3)}, \|\Pi_0\|_{L^{\infty}}$.

To derive the singularity for $\p_z \Pi$, we take the derive $\p_z$ for the equation by $\Pi$ and then obtain a Riccati type equation for $\p_z\Pi$
\be\label{a2}
\p_t \p_z\Pi+(u_r\p_r +u_z\p_z-2\Pi\p_z) \p_z \Pi- 2(\p_z\Pi)^2+ \p_z u_r\p_r \Pi+\p_z u_z\p_z \Pi=0.
\ee

Note that $u_r(t,r=0,z)\equiv 0$, we have $\p_z u_r(t,r=0,z)\equiv 0$. Hence if we restrict the equation (\ref{a2}) to $r=0$, then we obtain that
\be\label{a3}
\p_t \p_z\Pi(t,0,z)+(u_z-2\Pi)\p_z \p_z\Pi(t,0,z)-2(\p_z\Pi)^2(t,0,z)+(\p_z u_z \p_z \Pi)(t,0,z)=0.
\ee

By the divergent free condition, we have
\be\no
\p_z u_z(t,0,z)=-\lim_{r\rightarrow 0^+}(\p_r u_r+\f{1}{r} u_r)(t,r,z)=-2\p_r u_r(t,0,z).
\ee
Hence we obtain
\be\label{a4}
\p_t \p_z\Pi(t,0,z)+(u_z-2\Pi)\p_z \p_z\Pi(t,0,z)-2(\p_z\Pi)^2(t,0,z)-2(\p_r u_r \p_z \Pi)(t,0,z)=0.
\ee

Define the particle trajectory on the axis of symmetry $\phi(t,z)$ as follows
$$
\left\{ \aligned
&\f{d}{dt} \phi(t,z)= (u_z-2\Pi)(t,0,\phi(t,z)),\\
&\phi(0,z)=z.
\endaligned \right.
$$
Then by setting $f(t,z)=\p_z\Pi(t,0,\phi(t,z))$ and $g(t,z)=\p_r u_r(t,0,\phi(t,z))=\lim_{r\rightarrow 0}\f{u_r}{r}(t,0,\phi(t,z))$, we know that
\be\no
\f{d}{dt} f(t,z)&=& 2 f^2(t,z)- 2 g(t,z) f(t,z)\\\no
&\geq& f^2(t,z)-g^2(t,z).
\ee
Integrating over $[0,t]$, we have
\be\no
f(t,z)-f(0,z)&\geq& \int_0^t f^2(s,z) ds-\int_0^t g^2(s,z) ds\\\label{a15}
&\geq& \int_0^t f^2(s,z) ds-\int_0^t \|\f{u_r}{r}(s,\cdot)\|_{L^{\infty}}^2 ds.
\ee
Fix $z=0$ and set $y_0=f(0,0)=\p_z \Pi_0(0,0)$, then by employing the estimate (\ref{ur1}) in Lemma \ref{ur}, we obtain
\be\label{a16}
f(t,0)&\geq& \int_0^t f^2(s,0) ds+y_0-\int_0^t \left\|\f{u_r}{r}(s,\cdot)\right\|_{L^{\infty}}^2 ds\\\label{a17}
&\geq& \int_0^t f^2(s,0) ds+y_0-C_*(1+t)^{3/2} t^{1/2}.
\ee

Now take $y_0\geq10^4 C_*^2$ and $T_*=\f{4}{y_0}$, for $t\in [0,T_*]$, we have
\be\no
f(t,0)&\geq& \int_0^t f^2(s,0) ds+ y_0- 4 C_* \times \f{1}{100 C_*}\\\no
&\geq& \int_0^t f^2(s,0) ds+\f{1}{2} y_0.
\ee

Define a new function $F(t)=\int_0^t f^2(s,0) ds+\f{1}{2} y_0$, then $F(t)$ satisfies the following inequality
\be\no
&&F'(t)\geq F^2(t),\quad t\in [0,T_*],\\\no
&&F(0)=\f{1}{2} y_0.
\ee
Hence we have
\be\no
F(t)\geq \f{y_0}{2-t y_0},
\ee
which implies that
\be\no
\lim\sup_{t\rightarrow t_0} F(t)=\infty, \quad \lim\sup_{t\rightarrow t_0} f(t,0)=\lim\sup_{t\rightarrow t_0} \p_z \Pi(t,0,\phi(t,0))=\infty,
\ee
where $t_0=\f{2}{y_0}<T_*$.

Note that on the axis $r=0$, the equation for $\Pi$ can be reduced to
\be\no
\p_t \Pi(t,0,z)+ (u_z-2\Pi)(t,0,z)\p_z \Pi(t,0,z)=0.
\ee
By the definition of $\phi(t,z)$, we have $\f{d}{dt} \Pi(t,0,\phi(t,z))\equiv 0$. This enables us to get
\be\no
\Pi(t,0,\phi(t,0))= \Pi(0,0,\phi(0,0))=\Pi_0(0,0).
\ee
Hence, if we choose $\Pi_0(0,0)=J_0>0$, then
\be\no
\lim\sup_{t\rightarrow t_0} (\Pi \p_z \Pi)(t,0,\phi(t,0)) =\infty.
\ee
From the equation (\ref{a101}) for $\Omega$, we get
\be\label{a6}
2\Pi\p_z\Pi=-\p_t \Omega-(u_r\p_r +u_z\p_z) \Omega+(\p_r^2+\f{3}{r}\p_r+\p_z^2)\Omega.
\ee
Therefore we see that at least one of the terms on the right side in (\ref{a6}) blows up
\be\no
\lim\sup_{t\rightarrow t_0} \b(\p_t \Omega+(u_r\p_r+u_z\p_z)\Omega-(\p_r^2+\f{3}{r}\p_r+\p_z^2)\Omega\b)(t,0,\phi(t,0)) =\infty.
\ee
This contradicts to our assumption that (\ref{hmhd}) is globally well-posedness in some Sobolev space $H^m({\mathbb{R}^3})$ for $m>\f{7}{2}$. Hence the incompressible viscous Hall-MHD system without resistivity is not globally well-posedness in any Sobolev space $H^m({\mathbb{R}^3})$ for $m>\f{7}{2}$. $\square$

\subsection{Inviscid case $\nu=0$.}

{\it Proof of Theorem \ref{b}.} As in the proof of Theorem \ref{a}, we will argue by contradiction. Same argument as before shows that (\ref{a15}) also holds in the inviscid case, so
\be\label{a15}
f(t,z)-f(0,z)&\geq& \int_0^t f^2(s,z) ds-\int_0^t \left\|\f{u_r}{r}(s,\cdot)\right\|_{L^{\infty}}^2 ds.
\ee
Fix $z=0$ and set $y_0=f(0,0)=\p_z \Pi_0(0,0)$, then by employing the estimate (\ref{b9}), we obtain
\be\label{a16}
f(t,0)&\geq& \int_0^t f^2(s,0) ds+y_0-\int_0^t \left\|\f{u_r}{r}(s,\cdot)\right\|_{L^{\infty}}^2 ds\\\label{a17}
&\geq& \int_0^t f^2(s,0) ds+y_0-\tilde{C}_{*} t.
\ee

Now take $y_0\geq 4 \tilde{C}_{*}^{\f{1}{2}}$ and $\tilde{T}_{*}=\f{4}{y_0}\leq \tilde{C}_{*}^{-1/2}$, for $t\in [0,\tilde{T}_{*}]$, we have
\be\no
f(t,0)&\geq& \int_0^t f^2(s,0) ds+ y_0-  \tilde{C}_{*} \tilde{C}_{*}^{-1/2}\\\no
&\geq& \int_0^t f^2(s,0) ds+\f{1}{2} y_0.
\ee

Then $F(t)=\int_0^t f^2(s,0) ds+\f{1}{2} y_0$ satisfies the following inequality
\be\no
&&F'(t)\geq F^2(t),\quad t\in [0,\tilde{T}_{*}],\\\no
&&F(0)=\f{1}{2} y_0.
\ee
Hence we have
\be\no
F(t)\geq \f{y_0}{2-t y_0}.
\ee
This implies that
\be\no
\lim\sup_{t\rightarrow t_0} F(t)=\infty, \quad \lim\sup_{t\rightarrow t_0} f(t,0)=\lim\sup_{t\rightarrow t_0} \p_z \Pi(t,0,\phi(t,0))=\infty,
\ee
where $t_0=\f{2}{y_0}<\tilde{T}_{*}$.

As before, if we choose $\Pi_0(0,0)=J_0>0$, then
\be\no
\lim\sup_{t\rightarrow t_0} (\Pi \p_z \Pi)(t,0,\phi(t,0)) =\infty.
\ee
and also by (\ref{b1}), the velocity field will also blow up
\be\no
\lim\sup_{t\rightarrow t_0} \b(\p_t \Omega+(u_r\p_r+u_z\p_z)\Omega\b)(t,0,\phi(t,0)) =\infty.
\ee
This contradicts to our assumption that (\ref{hmhd}) is globally well-posedness in some Sobolev space $H^m({\mathbb{R}^3})$ for $m>\f{7}{2}$. Hence the incompressible viscous Hall-MHD system without resistivity is not globally well-posedness in any Sobolev space $H^m({\mathbb{R}^3})$ for $m>\f{7}{2}$. $\square$

\br\label{remark5}
{\it As one can see from the above proof, the convective term $(u_r\p_r+u_z\p_z)\Pi$ may prevent the shock formation. For the incompressible viscous Hall-MHD by restricting on the axis, we have good control on the gradient of $u_r$ and $u_z$, showing that the smoothing effect of the convective term is not strong enough and can not prevent the formation of shock-type singularity in the magnetic field. }
\er

{\bf Acknowledgements.} The second author would like to thank Prof. Zhouping Xin for the wonderful discussions.
D. Chae's research is supported partially by NRF Grants no.
 2006-0093854 and  no. 2009-0083521.

\bibliographystyle{plain}

\end{document}